# Cutting sequence and Sturmian sequence in billiard
# 台球中的切割序列和 Sturmian 序列

Zhiyu.Liu18

May 4, 2022

# Abstract

The winning rule of billiards is to drive the billiard ball on the table into the designated holes. We try to study the trajectory of the billiard ball, so that we can predict the direction of the ball. For rational slopes, we got cutting sequence by setting up the square torus. We simplified cutting sequence using shearing and flipping and we obtain the transformation between trajectory slope and cutting sequence. For irrational slopes, we look at some properties of Sturmian sequence, which help us distinguish between cutting sequence and Sturmian sequence. In conclusion, in the case of different slopes, we use different sequences to do research.

# 摘要

台球的获胜规则是把台球桌上的球打入指定的洞中。我们试图研究台球的运动轨迹，这样我们就可以预测球的运动方向。对于有理数斜率，通过建立方形环面得到切割顺序。我们利用剪切和翻转简化了切割顺序，得到了轨迹斜率与切割顺序之间的变换规则。对于无理数斜率，我们研究了 Sturmian 序列的一些性质，这有助于我们区分切割序列和 Sturmian 序列。综上所述，在不同斜率的情况下，我们使用不同的序列来进行研究。

# Contents

# 1. Introduction

Billiards are one of the most popular ball games in the world in recent years. In this sport, we hit the ball on the table with a long club in order to make the ball fall into the hole. According to Newton's first law of motion, in the absence of any external force, the ball will always move in a straight line. However, the ball does not usually move in a straight line on a square table and will change direction once it hits the edge of the table. What we're interested in is the billiard ball's trajectory. We can plot the trajectory of the ball step by step, using the angle of reflection is equal to the angle of incidence. Actually, drawing the full trajectory is tedious. Therefore, we try to figure out how to identify if the trajectory is periodic or non-periodic, and find the period. Actually, billiard tables are rectangular. To facilitate my research, we study with a square which is a special quadrilateral. We learned lots of useful methods and tools to express the trajectory in an easier way. We find when the slope of the trajectory is rational, the trajectory is periodic. First of all, we generated the torus table by square table to get a sequence named cutting sequence, and we used this periodic sequence to simplify complicated trajectories. Then we had a series of studies on cutting sequence, such as shear and flip. This led transform between trajectory slope and cutting sequence available. We only need to know one of them to derive the other. Moreover, the trajectory with irrational slope cannot be ignored. We get Sturmian sequence, and we learn its properties and compare it with cutting sequence we learned early.

# 2. A special billiard: square

First of all, Consider the most special case of a polygonal table: a square, because squares are centrosymmetric. To visualize the trajectory of the ball, let us assume that the ball is a little point, and there is no friction so that it runs forever. More importantly, when balls hit the edge of the table, the angle of reflection is equal to the angle of incidence. (We assume the ball never hits a vertex of square.) If we hit the ball vertically or horizontally, it will bounce between two points on parallel edges again and again. We say that this trajectory is periodic, because periodic trajectories mean that the ball repeats its path all the time. Once the ball starts moving, it goes back to the starting point after two collisions, so this trajectory has a period of 2. If we change the Angle of the strike, we can also draw the trajectory of the periods which are more than 2.

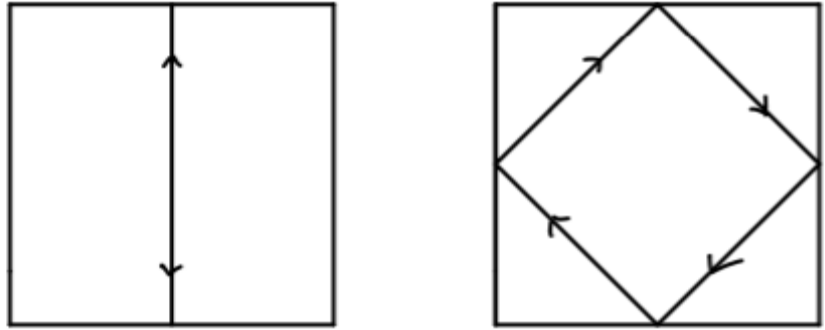

Figure 1: Trajectories with period 2 and 4

Can we find any trajectories with period 3? The answer is no. We hit a ball at an Angle of θ, then we can get ∠1=180°-2θ, ∠3=90°-θ, ∠4=180°-2∠3=2θ. Finally, we get ∠1+∠4=180°, and we know that the angles of a triangle add up to 180 degrees so there is no triangular in the trajectory. Similarly, we cannot find any trajectories with odd period.

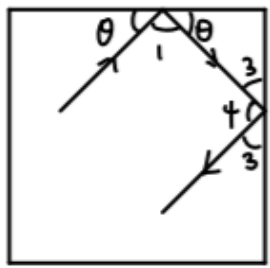

Figure 2

# 3. Unfolding the square table

However, as the number of periods increase, the trajectories we draw will get denser. It is impossible to study a non-periodic trajectory by drawing paths in a table. We can draw the trajectory in numerous tables, by unfolding the table.

## 3.1 Unfolding a trajectory into a straight line

Consider an example in Figure 3 below, and the trajectory of the part can be expressed as A→B→C→D. When the ball hits the top edge, we fold the table upward and draw a line segment AB' that is symmetric to the line segment AB. Then we continue to draw the trajectory of the ball inside the folded table. The ball is going to hit the right edge of the table so we get line segment B'C' for the right edge symmetry. By doing this over and over again, we can turn all the trajectories into a straight line. In the folding process, the trajectory itself invariant, we just change the position of the trajectory segments in different directions by symmetry. After four folds, we find that the ball is going along D' →A' which is same as it started from D→A. This means that the trajectory of the ball after it arrives A' will repeat the trajectory of the ball from A. In conclusion, the ball goes from A to A', completes a period of 4. [2]

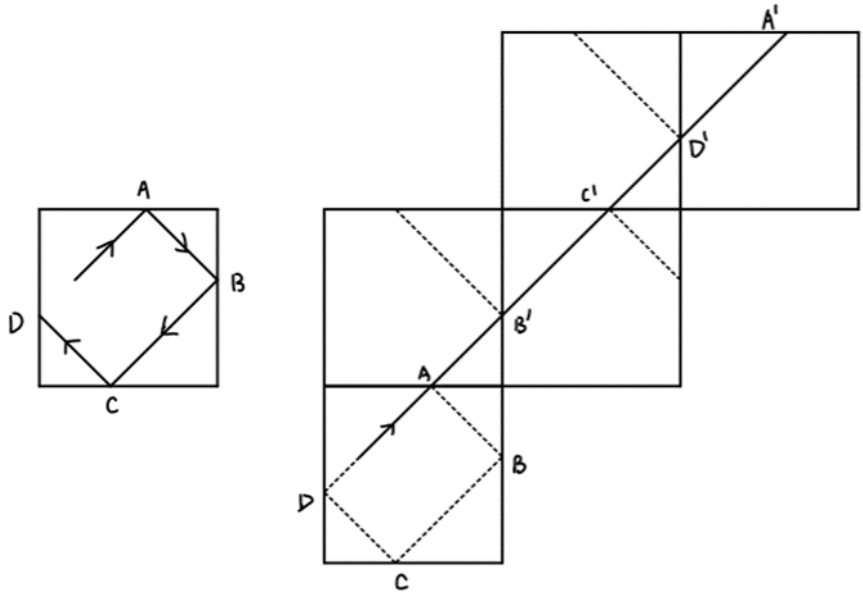

Figure 3

In the example, the trajectory is periodic. If I change the initial Angle(slope) of motion of the ball, will it still be periodic? No: If the trajectory is periodic, the ball must return to its original position after a period of movement. We try to calculate the slope in the example before (figure 4), which is equal to calculate ∠A'AO. It can be proved that the quadrilateral AOA'O' is a parallelogram because AO and A'O' are parallel and equal. Therefore, the slope is ∠A'AO=∠O'OM=$\frac{2}{2}$=1 = $\frac{P}{Q}$, where P is the number of times the square is unfolded upwards, and Q is the number of times the square is unfolded to the right. Moreover, P and Q must be rational numbers, so that the ratio of two rational numbers also must be a rational number. We say that a trajectory on the square billiard table is periodic if and only if its slope is rational.

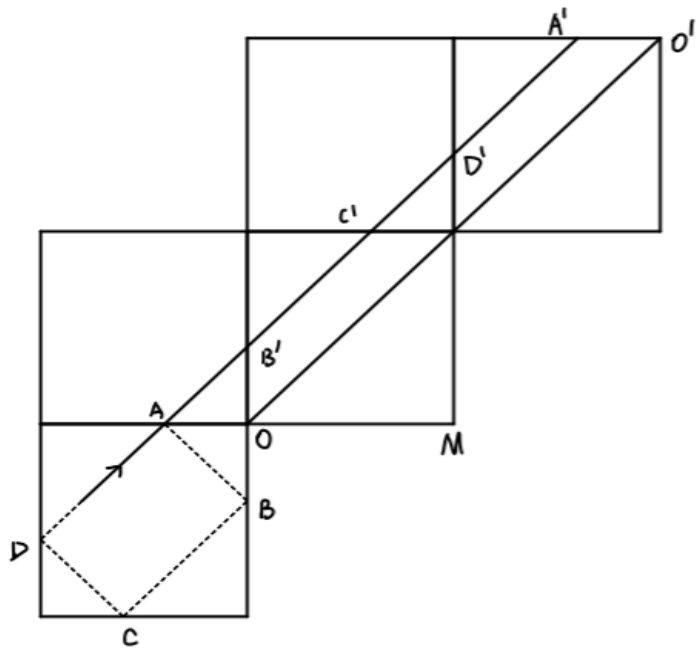

Figure 4

## 3.2 the square torus

If we fold a square upward, and then fold the resulting image to the right, we get a big square which is made up of four original squares. In this way, because of symmetry, the two edges above and below the new square all represent the edge A in the original square. When the ball hits the top edge of the new square, it reappears in the same place on the bottom edge. The left and right edges of a square have the same properties. Therefore, we can connect the top and bottom edges first and then the left and right edges. We finally get a torus, shown in figure 5. [2]

To understand it, thinking of the trajectory on a square torus as a three-dimensional trajectory on a surface, we can imagine a bug walking in a straight line on both planes. If the bug's journey on the square torus is horizontal, the 3D torus's equivalent path will be looped up like an equator and return to where it began (path A in Figure 5b). The corresponding path on the 3D torus will go through the hole and return to where the bug started if the bug's path on the square torus is vertical (path B in Figure 5b). If the slope of the path on the square torus is rational, the bug will loop around and around the torus and return to its starting position.

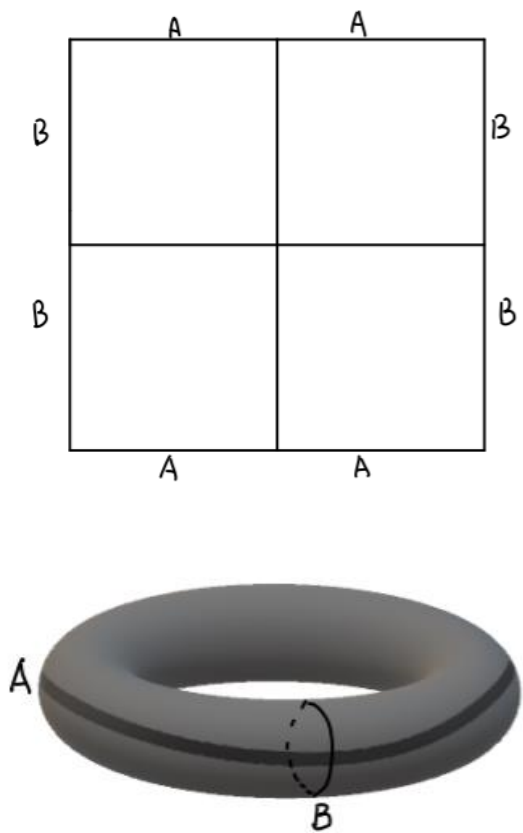

Figure 5: (a) the square torus (b) the square torus into a 3D torus

## 4. Cutting sequences

By construct the square torus, we represent the entire trajectory easier in a square with fewer line segments. For example, there are two lines in the square torus for trajectory with period 4(figure 6) and three lines for the trajectory with period 6.

How can we generate a sequence with respect to a trajectory on the square torus? The steps are followed: First, we choose a start point. Then, when the trajectory crosses the top or bottom edge, we record an A, and when it crosses the left or right edge, we record a B. finally, we obtain an infinite sequence named cutting sequence. If the trajectory is periodic, the corresponding cutting sequence is also periodic. For example, the trajectory in Figure 6 has cutting sequence $\overline{AB}$. We start moving at the bottom intersection with edge A, we record an A. Next, it crosses the right edge B, and the cutting sequence becomes AB. Then it reappears at the left edge B and crosses the top edge A. finally, the trajectory come back to the start. Therefore, the cutting sequence is …ABABAB…, it can be written as form $\overline{AB}$. The cutting sequence is not unique, it depends on the point where we started recording. If we started somewhere else, we would get another cutting sequence $\overline{BA}$. Actually, they represent the same trajectory.

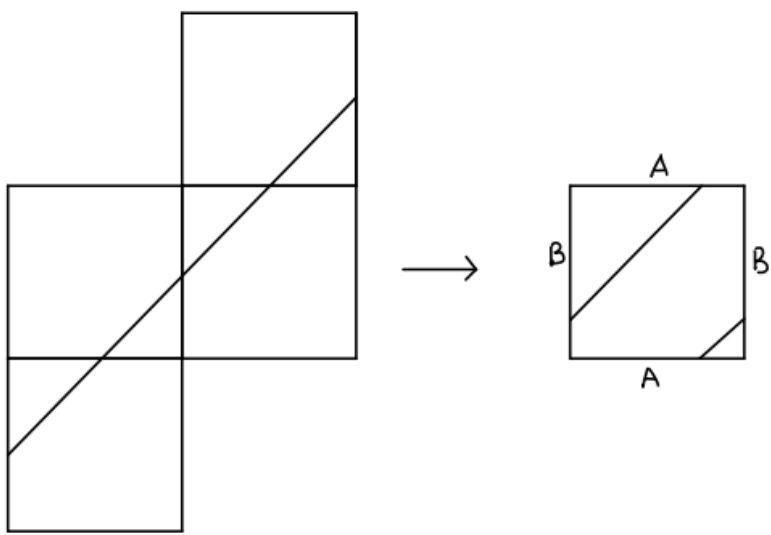

Figure 6: A trajectory ABAB with cutting sequence $\overline{AB}$

**Proposition 4.1.**[1] Consider a periodic trajectory on the square torus, and a billiard path in the same direction on the square billiard table. If the cutting sequence corresponding to the trajectory is $\overline{w}$, then the sequence of edges that the billiard ball hits in one period is ww.

In the length of one period of the torus trajectory, the billiard path traverses the same sequence of edges, twice. After hitting the sequence of edges w, the billiard path's orientation is reversed, thus it takes two cycles ww to return to the starting point in the same direction and complete a period.

Without considering the square tours, there is another way to generate a cutting sequence. Consider a trajectory in the square table. When line intersects a vertical side, we record B. When the line intersects a horizontal side, we record A. This way is quite similar to how we get Sturmian sequence in chapter 8.

Combined with the method we used in Section 3, the number of A in the cutting sequence is equal to the number of times the square is unfolded upwards, and the number of B in the cutting sequence is equal to the number of times the square is unfolded to the right. So, we can generate a new formula for the slope $\frac{p}{q}$ (in lowest terms), where the number of A is p and the number of B is q. The cutting sequence corresponding to the torus trajectory has period n=p+q. It is important to know the cutting sequence is not equal to the billiard path on the square table, and it is a simple expression of trajectory on the square torus. We can get the corresponding trajectory of the ball from cutting sequence and vice versa. Actually, the square torus consists of four squares. Therefore, the square torus A trajectory on the square torus can transform into a billiard path on the square table by folding the torus vertically and horizontally, and the period will double (figure 7). The billiard path has period 2n. This also shows that there are no odd periods in billiard.

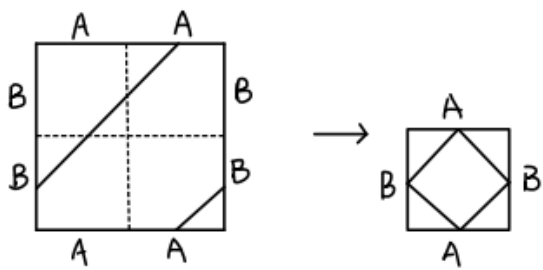

Figure 7: transformation between the square torus and the square table

Giving an arbitrary sequence, how can we check it is a cutting sequence? The following proposition and corollary could be used.

**Proposition 4.2.** [1] If an infinite sequence of As and Bs has two As in a row somewhere and also has two Bs in a row somewhere, then it is not a cutting sequence corresponding to a trajectory on the square torus.

Simple proof: We assume there exist a cutting sequence $\overline{BBAA}$, and try to draw it in a square torus. We find it is impossible, because the trajectories drawn intersect instead of parallel to each other.

**Corollary 4.3.** [1] A given cutting sequence on the square torus has blocks of multiple As separated by single Bs, or blocks of multiple Bs separated by single As, but not both.

  Suppose there is a sequence …AAABAAA…BBBABBB…. It is not a cutting sequence by proposition 4.2, because this sequence has both AA and BB.

## 5. shearing the square torus

  If we take a square torus and we cut it vertically, and then we twist it and dip the cut, we still get a square torus finally. Because of this property of the ring, we have the following theorem:

**Theorem 5.1.** [1] Given a trajectory $\tau$ on the square torus with slope greater than 1, and its corresponding cutting sequence c($\tau$), let $\tau$' be result of applying $\begin{bmatrix} 1 & 0 \\ -1 & 1 \end{bmatrix}$ to $\tau$. To obtain c($\tau$') from c($\tau$), shorten each string of As by 1.

  Proof by an example: Suppose there is a cutting sequence $\overline{BABAA}$(figure 8). First, we shear it by $\begin{bmatrix} 1 & 0 \\ -1 & 1 \end{bmatrix}$, so that make the square become a parallelogram. It does not affect edge B, and the slope of the trajectory and the A side has decreased by 1. We will explain how different shear work on the square torus in the later section. We call the diagonal of the parallelogram a, and the cutting sequence can be rewritten as $\overline{BABAaA}$. Then we cut the parallelogram diagonally, and put the two pieces together to form a new square torus. There is no A anymore, so we get $\overline{BBa}$. The top and bottom edges are actually the original diagonal a, and if we still use A to represent the top and bottom, we end up with $\overline{BBA}$. In conclusion, every part that only has A in it is going to subtract an A.

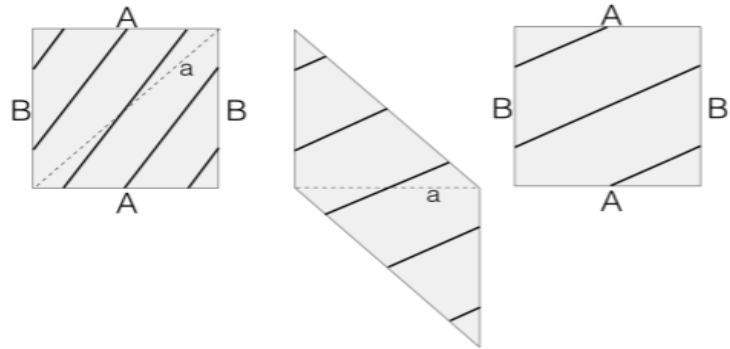

Figure 8: (a) The trajectory with cutting sequence $\overline{BABAA}$ ( (b) after shearing by $\begin{bmatrix} 1 & 0 \\ -1 & 1 \end{bmatrix}$ (c) reassembled back into a square with cutting sequence $\overline{BBA}$.[1]

In summary, we have the blew theorem which are practicable for all possible cutting sequences.

**Theorem 5.2.**[1] Given an infinite sequence of As and Bs, iterate the following process:

1. If it has AA somewhere and also BB somewhere, reject it; it is not a valid cutting sequence.

2. If it is …BBBABBB… or …AAABAAA…, reject it; it is not a valid cutting sequence.

3. If it has multiple As separated by single Bs, delete an A from each block.

4. If it has multiple Bs separated by single As, reverse As and Bs.

The above process can check if the given sequence is cutting sequence. When we keep doing these steps for a cutting sequence, we will never be rejected. Because when we shear with a cutting sequence, we still get a cutting sequence. But if we enter a valid sequence, the result is eventually rejected.

By shearing we can shorten the complex cutting sequence. However, when A is reduced by only one through the above steps, the cutting sequence still contains a large amount of B. We can switch A and B, then we can use the above Theorem 5.1 to simplify the cutting sequence.

Detail proof and more related materials for Theorem 5.1 are in [3], [4], [5] and [6].

## 6. continued fractions and cutting sequences

The continued fraction expansion gives an expanded expression of any given number, which can be used to apply the transformation of cutting sequence.

**Algorithm 6.1.**[1] We begin with the entire number as our "remainder," and iterate the following procedure.

1. If the remainder is more than 1, subtract 1.

2. If the remainder is between 0 and 1, take the reciprocal.

3. If the remainder is 0, stop.

We name the new expression of the number by a continued fraction.

**Example 6.2.** we write the continued fraction expansion of 7/4.

$$\frac{7}{4}=1+\frac{3}{4}=1+\frac{1}{4/3}=1+\frac{1}{1+\frac{1}{3}}$$

If we have a trajectory of a given rational slope, can we find the corresponding cutting sequence? Similarly, if can we find the slope of a trajectory form its cutting sequence. The answer is yes. We can do all this by shearing and flipping.

**Cutting sequence to continued fraction expansion:**

**Algorithm 6.3.**[1] Given a cutting sequence, perform the following procedure:

1. If the sequence has multiple As separated by single Bs, decrease the length of each string of As by 1. (We call it shear for short)

2. If the sequence has multiple Bs separated by single As, change all the Bs to As and all the As to Bs. (We call it flip for short)

3. If the sequence is an infinite string of Bs, stop. The slope is zero.

For 6.3.1, we apply the shear $\begin{bmatrix} 1 & 0 \\ -1 & 1 \end{bmatrix}$ to the square torus. Suppose there exist a vector $\begin{bmatrix} a \\ b \end{bmatrix}$. After applying, the vector becomes $\begin{bmatrix} a \\ b \end{bmatrix} * \begin{bmatrix} 1 & 0 \\ -1 & 1 \end{bmatrix} = \begin{bmatrix} a-b \\ b \end{bmatrix}$. The slope is $\frac{a-b}{b} = \frac{a}{b} - 1$, so the slope decreases 1.

For 6.3.2, we apply the flip $\begin{bmatrix} 0 & 1 \\ 1 & 0 \end{bmatrix}$, there is $\begin{bmatrix} a \\ b \end{bmatrix} * \begin{bmatrix} 0 & 1 \\ 1 & 0 \end{bmatrix} = \begin{bmatrix} b \\ a \end{bmatrix}$, and the slop becomes $\frac{b}{a}$, so the slope inverts.

In summary we have:

1. The shear subtracts 1 from the slope.

2. The flip inverts the slope.

Moreover, the trajectory corresponding to an infinite string of Bs has a slope of 0. Then we can find the slope, and we will do it in the following example.

**Example 6.4.** we perform algorithm 6.3 in cutting sequence $\overline{BABAABAABAA}$. Consider with algorithm 6.3, we do the following steps

shear → flip →shear→flip→shear→shear→shear

the cutting sequence changes as follows

$\overline{BABAABAABAA}$ —shear→ $\overline{BBABABA}$ —flip→ $\overline{AABABAB}$ —shear→ $\overline{ABBB}$

—flip→ $\overline{BAAA}$ —shear→ $\overline{BAA}$ —shear→ $\overline{BA}$ —shear→ $\overline{B}$

We solve it the other way around, because we know $\overline{B}$ has a slope of 0.

$$\frac{7}{4} = \cfrac{1}{\cfrac{1}{0+1+1+1}+1}+1$$

Finally, we find the trajectory corresponding to cutting sequence $\overline{BABAABAABAA}$ has a slope of 7/4. This is the same as the continued fraction expansion of 7/4.

**Continued fraction expansion to cutting sequence:**

How can we construct the cutting sequence corresponding to a trajectory with slope 7/4? We can think of it as the reverse step of example6.4.

**Algorithm 6.5.**[1] Given a slope, perform the following procedure:

1. Start with a string of Bs.

2. If the slope increases 1, insert one As between each pair of Bs.

3. If the slope invert, change all the Bs to As and all the As to Bs.

4. Repeat steps 2 and 3 to get the desired slope and the corresponding cutting sequence

We start at 0 and keep adding one and inverting, and we will end up with every slope that we want. During this period, we make corresponding changes to the cutting sequence and finally get the target cutting sequence.

The other thing we can say is that a cutting sequence only depends on the slope, it doesn't depend on where the ball start. If we hit the ball at the different positions in the billiard, as long as the slope of the trajectory is the same, we get the same cutting sequence. And that explains why, even though we initially ignored the situation that the ball hit the vertex, we can still use cutting sequence to figure out all the trajectory of the ball. We just have to change the starting point of the ball, and we can get the cutting sequence in this case as well.

Given a cutting sequence, we can use algorithm 6.3 to get the slope. Because, every time we change cutting sequence, the slope changes, and the slope will go to 0. We can start at 0 and work backwards to get the slope we started with. Similarly, if we have a slope, I can start with $\bar{B}$. and work backwards to get the original cutting sequence. We use the continued fraction expansion to realize the conversion between slope and cutting sequence. The important thing is that in this process, we didn't plot the trajectory in square torus like before. This means that even if we can't draw the trajectory of the ball in a polygon, we can still use cutting sequence to study its motion.

# 7. Every shear can be expressed by basic shears

So far, we have only considered the effect of the shear $\begin{bmatrix} 1 & 0 \\ -1 & 1 \end{bmatrix}$ on the cutting sequence corresponding to a trajectory whose slope is greater than 1. Can we use any other shear? Yes, but its effect on the slope of the trajectory maybe harder to state. Let us study from two basic shears.

**Proposition 7.1.**[1] The effects of applying the basic shears $\begin{bmatrix} 1 & 1 \\ 0 & 1 \end{bmatrix}$ and $\begin{bmatrix} 1 & 0 \\ 1 & 1 \end{bmatrix}$ to a linear trajectory on the square torus, with respect to the effect on the associated cutting sequence, are:

(a) $\begin{bmatrix} 1 & 0 \\ 1 & 1 \end{bmatrix}$: Lengthen every string of As by 1.

(b) $\begin{bmatrix} 1 & 1 \\ 0 & 1 \end{bmatrix}$: Lengthen every string of Bs by 1

Note that a "string of As" or a "string of Bs" can have length 0, and we add A or B directly.

Proof. (a) we know $\begin{bmatrix} 1 & 0 \\ 1 & 1 \end{bmatrix}$ is the inverse matrix of $\begin{bmatrix} 1 & 0 \\ -1 & 1 \end{bmatrix}$, and shear $\begin{bmatrix} 1 & 0 \\ -1 & 1 \end{bmatrix}$ shorten each string of As by 1. Instead, every string of As lengthen by 1.

(b) when we do shear $\begin{bmatrix} 1 & 0 \\ 1 & 1 \end{bmatrix}$, edge B does not change, and edge A rotates 45 degrees clockwise. When we do shear $\begin{bmatrix} 1 & 1 \\ 0 & 1 \end{bmatrix}$, edge A does not change, and edge B rotates 45 degrees clockwise. Compared with two situations, the roles of A and B are reversed.

Then we can reduce every shear to a composition of these two:

**Proposition 7.2.**[1] Every 2 × 2 matrix with nonnegative integer entries and determinant 1 is a product of powers of the basic shears $\begin{bmatrix} 1 & 1 \\ 0 & 1 \end{bmatrix}$ and $\begin{bmatrix} 1 & 0 \\ 1 & 1 \end{bmatrix}$. we only consider matrices with determinant 1. The reason is that we can reassemble the sheared torus back into the original square, it will not change the area.

Simple proof: suppose there is a matrix $\begin{bmatrix} a & b \\ c & d \end{bmatrix}$ with determinant 1, which means that ad-cb=1. We find that

$$\begin{bmatrix} a & b \\ c & d \end{bmatrix} * \begin{bmatrix} 1 & 1 \\ 0 & 1 \end{bmatrix}^{-1} = \begin{bmatrix} a-c & b-d \\ c & d \end{bmatrix},$$

and the determinant is (a-c)*d-(b-d)*d=ad-cb=1. Also,

$$\begin{bmatrix} a & b \\ c & d \end{bmatrix} * \begin{bmatrix} 1 & 0 \\ 1 & 1 \end{bmatrix}^{-1} = \begin{bmatrix} a & b \\ -a+c & -b+d \end{bmatrix},$$

and the determinant is 1. Once we have a matrix, we can let it time $\begin{bmatrix} 1 & 1 \\ 0 & 1 \end{bmatrix}^{-1}$ to make its element second column smaller, and time $\begin{bmatrix} 1 & 0 \\ 1 & 1 \end{bmatrix}^{-1}$ to make its element in first second column smaller. Finally, we will get a matrix $\begin{bmatrix} 1 & 1 \\ 0 & 1 \end{bmatrix}$ or $\begin{bmatrix} 1 & 0 \\ 1 & 1 \end{bmatrix}$. Multiply both sides of this equation by some $\begin{bmatrix} 1 & 1 \\ 0 & 1 \end{bmatrix}$ or $\begin{bmatrix} 1 & 0 \\ 1 & 1 \end{bmatrix}$, and we get rid of these inverse matrixes.

**Example 7.3.** there is a trajectory of slope 1, with cutting sequence $\overline{AB}$ on the square torus. What is the new cutting sequence after shearing the torus via the matrix $\begin{bmatrix} 3 & 1 \\ 2 & 1 \end{bmatrix}$?

Firstly, we express the matrix $\begin{bmatrix} 3 & 1 \\ 2 & 1 \end{bmatrix} = \begin{bmatrix} 1 & 1 \\ 0 & 1 \end{bmatrix}\begin{bmatrix} 1 & 0 \\ 1 & 1 \end{bmatrix}\begin{bmatrix} 1 & 1 \\ 0 & 1 \end{bmatrix}$. So, the question becomes shearing the torus via the matrix $\begin{bmatrix} 1 & 1 \\ 0 & 1 \end{bmatrix}, \begin{bmatrix} 1 & 0 \\ 1 & 1 \end{bmatrix}, \begin{bmatrix} 1 & 1 \\ 0 & 1 \end{bmatrix}$ one by one. Reminding the Proposition 7.1., we have $\overline{AB} \to \overline{ABB} \to \overline{AABAB} \to \overline{ABABBABB}$.

We have seen how important shearing is in the study of cutting sequence so far. With the help of shearing, we relate the slope of the trajectory to cutting sequence, and we just need to know one of them to get the other one without plotting. However, our series of studies are periodic trajectories and the corresponding cutting sequence. How about non-periodic trajectories, let us study it in next chapter.

## 8. Sturmian sequence. Basic properties

Recall: cutting sequence in a square table are coding of trajectory of rational slope obtained by labeling by A and B respectively its horizontal and vertical sides. This is a periodic sequence which generated by A and B.

Similarly, the most intuitive method to get a Sturmian sequence is to consider a line in the square table with an irrational slope and construct a sequence by considering its intersections with an integer grid. When the line intersects a vertical side, we record 1. When the line intersects a horizontal side, we record 0. We have already proof that the trajectory of a ball struck from an irrational angle in billiard has no period. Therefore, the Sturmian sequence is infinite non-periodic sequence which consist of 0 and 1. Sturmian sequences are defined by a rich classical theory that ties them to continuing fraction expansions of the slopes of the trajectories [7]. Christoffel [8] and Smith [9] explored Sturmian sequences in the 1870s, and they are significant examples in Morse and Hedlund's theory of symbolic dynamics [10]. They have been characterized as the non-periodic sequences of minimal complexity [11].

**Example 8.1.** find the Sturmian sequence with slop $\sqrt{2}$.

We can use matlab to plot a graph:

*>>fplot(@(x)2^(0.5)*x+0.5,[0,50])*

*grid on*

*axis equal*

we get the graph below:

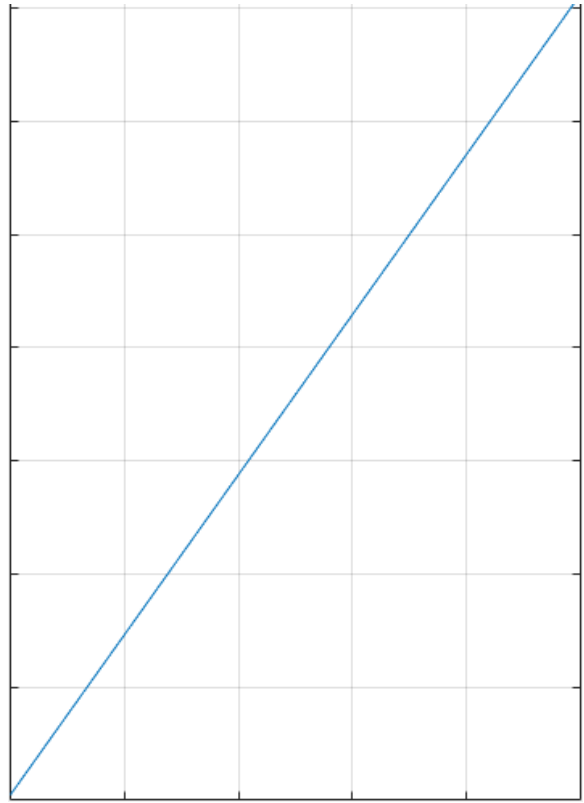

Figure 9: trajectory with irrational slope $\sqrt{2}$

We get a Sturmian sequence 101010010100…

**Definition 8.2.**[1] The complexity function $p(n)$ on a sequence is the number of different "words" of length n in the sequence.

There is an easy example to understand complexity. If we want to find p(2) for "door", we need to find all the different combinations with two letters. That is "do", "oo" and

"or". We get p(2)=3.

We know the highest possible complexity for a sequence of As and Bs is $p(n)=2^n$, because any element has only two choices (A or B). Periodic cutting sequence which has m elements in a cycle on the square torus have complexity p(n)=n+1 for n<m and complexity p(n)=m for n≥m.

**Example 8.3.** calculate complexity p(n) for the cutting sequence $\overline{ABBA}$.

For n=1, we have "A",'B". p(1)=2;

For n=2, we have "AB","BB","BA". P(2)=3;

For n=3, we have "ABB","BBA","BAA","AAB". P(3)=4;

For n=4, we have "ABBA","BBAA","BAAB","AABB". P(4)=4;

For n=5, we have "ABBAA,"BBAAB","BAABB","AABBA". P(5)=4.

For larger n, p(n) is always equal to 4. Why? When we choose n is bigger than 5, we collect the subsequence with length m one by one. Firstly, we can find a subsequence with length m starting with A, followed by B. When we get to the (m+1) th sequence, we find that this sequence is the same as the first sequence I got. In a word, we can only find at most m different "words" of length n in the cutting sequence for n≥m.

**Definition 8.4.**[12] A sequence u is called Sturmian if it has complexity $P_u(n) = n + 1$.

We can get some idea form example 8.3. Sturmian u is non-periodic, so we can consider this case as m is infinite. That means n<m always exists. Therefore, we get p(n)=n+1.

**Proposition 8.5.**[12] A Sturmian sequence is recurrent, that is, every word that occurs in the sequence occurs an infinite number of times.

Proof.   Suppose there is a word U with length n, occurs in a Sturmian sequence u a finite number of times. So there exist a subsequence v will not contain U. We know $p_u(n)$=n+1, and v is a part of u. Therefore $p_v(n)$≤n, which means v is periodic by example 8.1. This is a contradiction.

**Lemma 8.6.**[12] If u is Sturmian, then exactly one of the words 00, 11 does not occur in u.

Proof.    In chapter 5, we briefly show that this lemma is also true for cutting sequence. Actually, the principle is the same.

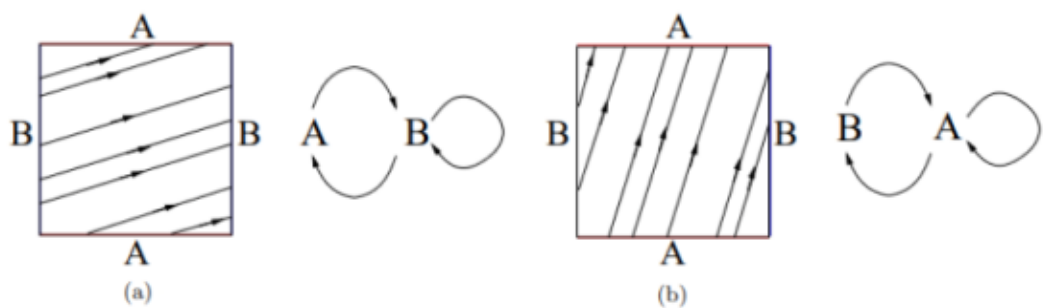

Figure 10: Possible transitions in the square torus

We can assume we have a trajectory in a rational angle θ with $0 \leq \theta \leq \pi/2$. We consider this question in two conditions: $0 \leq \theta \leq \pi/4$ and $\pi/4 \leq \theta \leq \pi/2$. If $0 \leq \theta \leq \pi/4$, as in figure 10(a), the subword AA does not occur in the cutting sequence. If there is a A in the cutting sequence, the next element must be B. And if there is a B in the cutting sequence, the next element can be A or B. Similarly, we see figure 10(b) for $\pi/4 \leq \theta \leq \pi/2$, it does not contain the subword BB. This case can be reduced to $0 \leq \theta \leq \pi/4$ by interchanging the role of A and B [7]. we can conclude that either there is only AA or there is only BB in the cutting sequence, the two cannot exist at the same time. For the irrational angle θ, these conclusions still hold. AA and BB is same as 00 and 11.

Also, we can proof it in another way. There are three different words of length 2 in Sturmian sequence u, because $p_u(2)=3$. The all combinations of 0 and 1 are 00,11,01 and 10. By the previous proposition, both 1 and 0 appear in u an infinite number of times, which means that 01 and 10 must occur in u. so exactly one of 00 and 11 must occur.

Further thinking: For $0 \leq \theta \leq \pi/4$, there is no BBB and AAA in the sequence. For $\pi/4 \leq \theta \leq \pi/2$, there is no AAA and BBB in the sequence.

**Definition 8.7.**[12] We say that a Sturmian sequence is of type 0 if 1 is isolated, that is, if 11 does not occur in the sequence. This is equivalent to saying that 0 occurs more frequently than 1. We say that the sequence is of type 1 if 00 does not occur.

By this definition, we divided the sequence into two categories. For example 8.1, it is a type 0 Sturmian sequence.

To better introduce the next lemma, we introduce a definition first.

Note: if U is a finite word, we denote by $|U|_a$ the number of occurrences of the letter a in U.

**Definition 8.8.** [12] A sequence u consists of 1 and 0 is balanced if, for any pair of words U, V of the same length occurring in u, we have $||U|_1 - |V|_1| \leq 1$.

Let us take some examples. For a periodic sequence $\overline{100}$, we list all different subword of different length.

When length is 1, we can choose U is 0 and V is 1(or U is 1 and V is 0), then $||U|_1 - |V|_1|=1$.

When length is 2, we have 3 choices 10, 00 and 01. We choose two of the three choices as U and V. therefore, $||U|_1 - |V|_1|=1$ or 0.

When length is 3, we have 3 choices 100,001,010. $||U|_1 - |V|_1|=0$.

When length is bigger than 3, we always have three choices, the value of $||U|_1 - |V|_1|$ is 0 or 1. When length is a multiple of 3, it can only be 0.

Cutting sequence is a balanced sequence.

**Lemma 8.9.** [12] If the sequence u is not balanced, there is a (possibly empty) word W such that 0W0 and 1W1 occur in u.

By definition 8.8, if the sequence u is not balanced, there are some pair of words U,V of the same length n occurring in u, we have $||U|_1 - |V|_1| >1$. That means we can find U,V which make $||U|_1 - |V|_1| =2$. There are two 1's in U while there are two 0's in V. Suppose now that A and B are words of minimal length with this property. Write

$A=a_0a_1\ldots a_{n-1}$ and $B=b_0b_1\ldots b_{n-1}$. By lemma 8.6, it cannot be 11 in A and 00 in B. If we let $a_0=a_{n-1}=1$, we must have $b_0=b_{n-1}=1$. If we have situation like $a_0=a_{n-1}=1$ and $b_0=b_1=0$, we can find a shorter pair by removing some prefix. It will not be minimal.

**Theorem 8.10.**[12] A sequence u is Sturmian if and only if it is a non-eventually periodic balanced sequence over two letters.

The complexity $p_u$ is an increasing function. If u is eventually periodic, then $p_u$ is bounded and there is an n such that $p_u(n + 1) = p_u(n)$. We know from previous studies that the cutting sequence is eventually periodic. If u is not eventually periodic, we should have $p_u(n) \geq n+1$, due to $p_u(1) \geq 2$ and $p_u$ is increasing. Otherwise, u will be a constant or an infinite sequence of identical elements.

For a Sturmian sequence, we have $p(n)=n+1$. The equation $P(n+1)=n+2>n+1>p(n)$ holds all the time. Thus, Sturmian sequence is non-eventually.

We can use lemma8.9 to check Sturmian sequence is a balanced sequence in a simple case. Suppose there is a Sturmian sequence $1W1\ldots 0W0$. It is not balanced, because we have $||1W1|_1 - |0W0|_1|=2$. We know there is only one of 11 and 00. Firstly, we choose the first element in W is 0. So, 00 exist, and there is no 11. The final element of W must be 0. Could W be 0? By the thinking under the lemma8.6, it is impossible. Maybe we can choose W=010. Then we get a Sturmian sequence $10101\ldots 00100$. However, this Sturmian sequence does not exist. We get a contradiction.

We draw the trajectory of 10101(Figure 11(a)) and 00100(Figure 11(b)) in the square table. Once 10101 appears in a Sturmian sequence, we get that the slope of the corresponding trajectory is less than 2/3. If there is a subsequence 00100 in a Sturmian sequence, the slope of the corresponding trajectory is greater than 2/3. The trajectory is a straight line, so the slope of the trajectory is constant. There is no constant that is both greater than two-thirds and less than two-thirds. That means 10101 and 00100 can not occur in a same Sturmian sequence at the same time. Therefore, we cannot find a non-balanced Sturmian sequence.

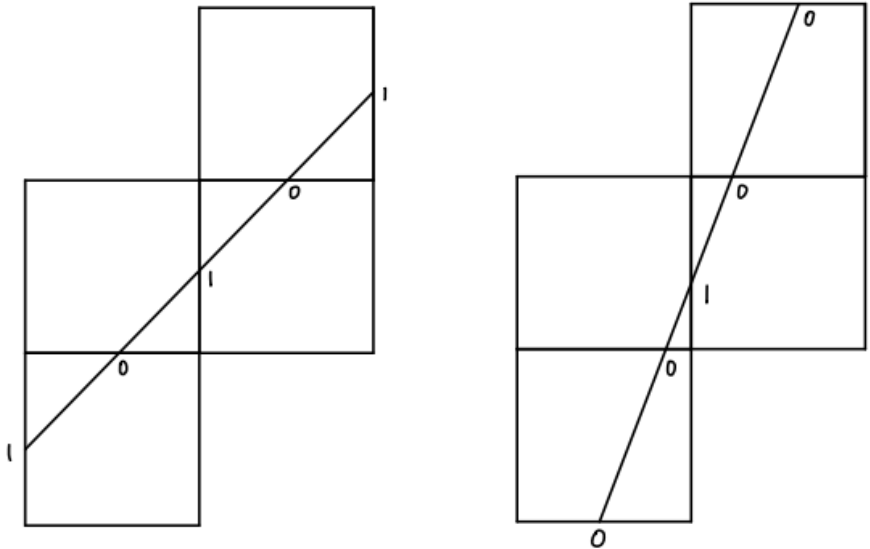

Figure 11:(a) trajectory of 10101 (b) trajectory of 00100

# 9.Continued fractions and Sturmian sequences

In chapter 6, we find continued fraction is useful in study for cutting sequence. We expand every rational number into continued fraction expansion. Can we do this with an arbitrary irrational number? The answer is yes.

**Example 9.1.** find the continued fraction expansion for $\sqrt{2}$.

We know that $1<\sqrt{2}<2$, so we suppose that

$$\sqrt{2}=1+x$$

Then square both sides of this equation

$$2=(1+x)2=1+x2+2x=1+x(x+2)$$

We can express x

$$x=\frac{1}{2+x}$$

So, we plug in $x=\frac{1}{2+x}$ on the right-hand side

We get

$$x = \cfrac{1}{2 + \cfrac{1}{2 + \cfrac{1}{2 + \cfrac{1}{2 + \cdots}}}}$$

Therefore,

$$\sqrt{2} = 1 + \cfrac{1}{2 + \cfrac{1}{2 + \cfrac{1}{2 + \cfrac{1}{2 + \cdots}}}}$$

We can't write it out in its entirety because its length is infinite. And that corresponds the Sturmian sequence with an irrational slope is non-periodic. That means the ball in the billiard will keep running all the time unpredictably.

Moreover, can we translate between a Sturmian sequence, and the slope of the corresponding trajectory? Shearing and flipping are also available for Sturmian sequence.

If there is a Sturmian sequence, it might be much difficult to know the slope. When we doing the translation with cutting sequence, we know the total sequence, and we can decrease it step by step. However, for a Sturmian sequence, the algorithm continues forever.

If we have a slope, it is also hard to generate the Sturmian sequence. Reminder: when we carry out the algorithm 6.3, we need an original sequence which has corresponding slope $a_k$. Because for every rational number, we have slope $= a_1 + \cfrac{1}{\cdots + \cfrac{1}{a_k}}$, for irrational number we will never find this $a_k$.

# 10. Compare with cutting sequence and Sturmian sequences

When we put Cutting sequence and Sturmian sequence together, we will find that they have many similarities and differences.

Let us do similarities first. They can represent the trajectory of the ball on the square table very well, and they all follow the laws of billiard ball movement. Also, they are biinfinite sequence which means they can go back and forth indefinitely.

The most obvious difference between the two is that cutting sequence is periodic and Sturmian sequence is non-periodic. The former is the trajectory of the ball from a rational angle, the latter is the trajectory of the ball from an irrational angle.

One way to identify them is to check the complexity function p(n). if p(n)=n+1 for all n, this sequence is a Sturmian sequence. If p(n+1) =p(n) happens, it is a cutting sequence.

## 11.conclusion

We start with the billiard ball in reality, and want to explore the mathematical truth contained therein. It is important to ignore physical factors like friction to consider the full trajectory of a ball on a billiard table. So, we began our study of the most special quadrilateral square. According to the theorem that the Angle of reflection is equal to the Angle of incidence, we can draw the trajectory of the ball in a square table. When the period is very small, it is easy to plot all the trajectories. And then when the period gets bigger, and even when the trajectory doesn't have a period, all of this becomes difficult. We were trying to represent all the trajectories of the ball in a more intuitive way. When the trajectory is periodic, the slope of the trajectory is rational. We constructed the square torus, and we got cutting sequence with two element A and B. We represent the path of the billiard ball as a sequence, so that we can do more mathematically related research. By shearing and flipping, we can make an efficient change to cutting sequence. And with the continued fractions expansion, we were able to relate the slope of the trajectory to cutting sequence. We do not have to draw the trajectories anymore, just know the slope and we will get cutting sequence and then all the trajectories are obvious. Furthermore, the trajectory of the ball may appear to have no period, and we can also use Sturmian sequence to represent its trajectory. After studying the properties of Sturmian sequence and cutting sequence, there are many similarities and differences. Their related properties can help us better understand the trajectory of the billiard ball. We can reject valid sequences by complexity function or theorems about balanced and non-eventually periodic.